\documentclass[12pt]{amsart}
\input{amssymb.sty}
\title{Boundaries of Reduced Free Group $C^*$-Algebras}
\author{Narutaka OZAWA}
\address{Department of Mathematical Sciences,
University of Tokyo, Komaba, 153-8914\\ \indent 
Department of Mathematics, UCLA, Los Angeles, CA 90095-1555}
\email{narutaka@ms.u-tokyo.ac.jp}
\date{November 20, 2004}
\thanks{The author was supported in part by JSPS}
\subjclass{Primary 46L05; Secondary 20F67}

\keywords{free groups, reduced group $C^*$-algebras, 
injective envelopes}
\newtheorem{thm}{Theorem}
\newtheorem{prop}[thm]{Proposition}
\newtheorem{cor}[thm]{Corollary}

\newcommand{\G}{\Gamma}
\newcommand{\bG}{\partial\Gamma}
\newcommand{\cg}{C^*_{\mathrm{r}}\G}
\newcommand{\cgb}{C^*_{\mathrm{r}}(\G,\bG)}
\newcommand{\vng}{\G\ltimes L^\infty(\bG,\mu)}
\newcommand{\p}{\varphi}

\newcommand{\C}{{\mathbb C}}
\newcommand{\M}{\mathcal{M}}
\newcommand{\id}{\mathrm{id}}
\begin{document}
\begin{abstract}
We prove that the crossed product $C^*$-algebra 
$C^*_{\mathrm{r}}(\Gamma,\partial\Gamma)$ of a free group $\G$ 
with its boundary $\partial\Gamma$ naturally sits 
between the reduced group $C^*$-algebra $C^*_{\mathrm{r}}\Gamma$ 
and its injective envelope $I(C^*_{\mathrm{r}}\Gamma)$. 
In other words, we have natural inclusion 
$C^*_{\mathrm{r}}\Gamma\subset C^*_{\mathrm{r}}(\Gamma,\partial\Gamma)\subset 
I(C^*_{\mathrm{r}}\Gamma)$ of $C^*$-algebras.
\end{abstract}
\maketitle
\section{Introduction}
A deep theorem of Kirchberg \cite{kirchberg}\cite{kp} 
states that for any separable exact $C^*$-algebra $B$, 
there exists a separable nuclear $C^*$-algebra $A$
that contains $B$ as a $C^*$-subalgebra. 
(Indeed, one can take $A$ as the Cuntz algebra $\mathcal{O}_2$.) 
In this note, we are interested in finding 
a more tight (or canonical) embedding of a given separable exact 
$C^*$-algebra into a separable nuclear $C^*$-algebra. 
This is partly motivated by Hamana's theory of injective envelopes. 
An embedding of $B$ into $A$ is called rigid if 
the identity map on $A$ is the only completely positive map on $A$ 
which restricts to the identity map on $B$. 
Hamana \cite{hamana} proved that for any $C^*$-algebra $B$, 
there exists an injective $C^*$-algebra $I(B)$, 
called the injective envelope of $B$, 
with a rigid embedding of $B$ into $I(B)$.  
We note that the injective envelope $I(B)$ is unique up to 
natural identification. 

Let $\G$ be the free group of rank $2\le r<\infty$ 
and $\bG$ be its boundary, i.e., 
the compact topological space of all one-sided 
infinite reduced words, equipped with the product topology. 
The free group $\G$ acts on $\bG$ by translation and 
this action is topologically amenable (cf.\ \cite{ar}). 
The reduced group $C^*$-algebra $\cg$ is naturally 
contained in the reduced crossed product $C^*$-algebra $\cgb$,
which is nuclear. 
A measure $\mu$ on $\bG$ is called \emph{quasi-invariant} 
if for any measurable subset $X\subset\bG$ and any $s\in\G$, 
one has $\mu(sX)=0$ if and only if $\mu(X)=0$. 
A measure $\mu$ on $\bG$ is called \emph{doubly-ergodic} 
if the diagonal action of $\G$ on 
$(\bG^2,\mu^{\otimes2})$ is ergodic. 
Fix a quasi-invariant and doubly-ergodic measure 
$\mu$ on $\bG$. 
For instance, the harmonic measure of any 
non-degenerate symmetric random walk on $\G$ 
has this property (Theorem 2.8 in \cite{kaimanovich}).
Associated with $\mu$, there is a natural 
inclusion of $\cgb$ into the crossed product 
von Neumann algebra $\vng$, which is injective; 
\[
\cg\subset\cgb\subset\vng.
\]
See \cite{okayasu} \cite{rr} for information on 
the type of the factor $\vng$.

We state the main theorem of this note. 
A somewhat similar result was proved in \cite{spielberg}.
\begin{thm}\label{thm}
Let $\G$ be the free group of rank $2\le r<\infty$, 
$\bG$ be its boundary and $\mu$ be a quasi-invariant 
and doubly-ergodic measure on $\bG$. 
If 
\[
\theta\colon\cgb\to\vng
\]
is a completely positive map with $\theta|_{\cg}=\id_{\cg}$, 
then $\theta=\id$. 
\end{thm}

\begin{cor}\label{cor}
Let $\G$ and $\bG$ be as above. 
Then, the nuclear $C^*$-algebra $\cgb$ sits (as a $C^*$-algebra) 
between $\cg$ and its injective envelope $I(\cg)$, 
i.e., we have natural inclusions 
$\cg\subset\cgb\subset I(\cg)$ of C$^*$-algebras.
\end{cor}
The author conjectures that for 
any separable exact $C^*$-algebra $B$, 
there exists a nuclear $C^*$-algebra between $B$ and 
its injective envelope $I(B)$. 

The theorem is a consequence of the following 
proposition from ergodic theory, inspired by \cite{ms}.
We recall that a map between spaces with $\G$-actions 
is called $\G$-equivariant if it commutes with 
the $\G$-actions. 
\begin{prop}\label{prop}
Let $\G$ be the free group of rank $2\le r<\infty$, 
$\bG$ be its boundary and $\mu$ be a quasi-invariant 
and doubly-ergodic measure on $\bG$. 
If 
\[
\p\colon C(\bG)\to L^\infty(\bG,\mu)
\] 
is a unital positive $\G$-equivariant map, 
then $\p=\id$. 
\end{prop}
Actually the same result \textit{mutatis mutandis} holds 
for a discrete subgroup $\G$ of the group of isometries on 
an exponentially-bounded Gromov hyperbolic space $\mathcal{X}$ 
and a closed subspace of $\partial\mathcal{X}$ 
on which $\G$ acts minimally. 
It is plausible that the assumption on $\mu$ can be relaxed 
to that $\mu$ is quasi-invariant and diffuse. 
(For this matter, it suffices to show this assertion under 
the additional assumption that $\p$ is a $\G$-equivariant 
$*$-homomorphism.)
\section{Proof of Proposition}
Let $\G$ be the free group of rank $2\le r<\infty$, 
$\bG$ be its boundary and $\mu$ be a quasi-invariant 
and doubly-ergodic measure on $\bG$. 
Let a unital positive $\G$-equivariant map
$\p\colon C(\bG)\to L^\infty(\bG,\mu)$ be given. 

We denote by $\M(\bG)$ the space of all finite Borel measures on $\bG$ 
and let $\M_{\le 2}(\bG)$ be its subset of measures 
which are supported on at most two points.
It is well-known (cf.\ Section 5 in \cite{adams}) that 
$\M_{\geq3}(\bG)=\M(\bG)\setminus\M_{\le 2}(\bG)$ is a tame $\G$-space. 
From this and the double ergodicity of $\mu$, it follows that 
any $\G$-equivariant Borel map from $\bG^2$ into $\M_{\geq3}(\bG)$ 
has image in a single orbit $\mu^{\otimes2}$-a.e. 
Since every non-trivial element of $\G$ acts hyperbolically on $\bG$, 
the stabilizer of this orbit is trivial. 
Since $\mu^{\otimes2}$ is diffuse, the pre-image of any point of 
this orbit can be divided into disjoint non-null pieces, contradicting 
ergodicity of $\mu^{\otimes2}$. 
Thus any $\G$-equivariant Borel map from $\bG^2$ into $\M(\bG)$ 
has image in $\M_{\le 2}(\bG)$ $\mu^{\otimes2}$-a.e. 

Fix a dense $\G$-invariant $*$-subalgebra $\mathcal{C}$ 
in $C(\bG)$ which is algebraically (over $\C$) generated 
by a countable set. 
Then, there exists a $\G$-equivariant Borel map 
\[
\p_*\colon\bG\ni\xi\mapsto\p_*^\xi\in\M(\bG)
\]
such that for $\mu$-a.e.\ $\xi\in\bG$, we have
\[
\forall f\in\mathcal{C}\quad
\int f(\eta)\, d\p_*^\xi(\eta) = \p(f)(\xi).
\]
We consider the $\G$-equivariant Borel map 
\[
\bG^2\ni (\xi,\eta)\mapsto
\p_*^\xi+\delta_\xi+\delta_\eta\in\M(\bG),
\]
where $\delta_\xi$ is the Dirac measure on $\xi\in\bG$. 
By the previous result, this map 
has image in $\M_{\le 2}(\bG)$ $\mu^{\otimes2}$-a.e. 
Since the diagonal subset of $\bG^2$ is $\mu^{\otimes2}$-null 
and $\xi$ and $\eta$ are independent, 
we conclude that $\p_*^\xi=\delta_\xi$ 
for $\mu$-a.e.\ $\xi\in\bG$. 
This means that $\p=\id$. 
\hspace*{\fill}$\Box$
\section{Proof of Theorem}
Let $\G$ be the free group of 
rank $2\le r<\infty$, $\bG$ be its boundary and $\mu$ be a 
quasi-invariant and doubly-ergodic measure on $\bG$. 
For $s\in\G$ and $f\in L^\infty(\bG,\mu)$, 
we define $s\cdot f\in L^\infty(\bG,\mu)$
by $(s\cdot f)(\xi)=f(s^{-1}\xi)$. 
The crossed product von Neumann algebra $\vng$ 
is generated (as a von Neumann algebra) by 
a unitary representation $\lambda$ of $\G$ and $L^\infty(\bG,\mu)$. 
They satisfy the covariance property that 
$\lambda(s)f\lambda(s)^*=s\cdot f$ 
for $s\in\G$ and $f\in L^\infty(\bG,\mu)$. 
Moreover, there exists a normal faithful 
conditional expectation $E$ from $\vng$ onto $L^\infty(\bG,\mu)$ 
such that for $s\in\G$ and $f\in L^\infty(\bG,\mu)$, one has 
$E(\lambda(s)f)=f$ if $s=1$ and $E(\lambda(s)f)=0$ if $s\neq1$. 
The $C^*$-subalgebra in $\vng$ generated by 
$\lambda(\G)$ and $C(\bG)$ is naturally $*$-isomorphic 
to the reduced crossed product $C^*$-algebra $\cgb$ 
and the $C^*$-subalgebra in $\cgb$ generated by 
$\lambda(\G)$ is naturally $*$-isomorphic 
to the reduced group $C^*$-algebra $\cg$; 
\[
\cg\subset\cgb\subset\vng.
\]

We prove Theorem \ref{thm} and Corollary \ref{cor}. 
\begin{proof}[Proof of Theorem \ref{thm}]
Let $\theta\colon\cgb\to\vng$ be a completely positive map 
such that $\theta|_{\cg}=\id_{\cg}$. 
Consider the unital positive map $\p=E\circ\theta|_{C(\bG)}$. 
Since $\theta|_{\cg}=\id_{\cg}$, the completely positive map 
$\theta$ is automatically a $\cg$-bimodule map (cf.\ \cite{ce}). 
Hence, for any $s\in\G$ and $f\in C(\bG)$, we have
\begin{align*}
\p(s\cdot f) &= E(\theta(\lambda(s)f\lambda(s)^*))\\
 &= E(\lambda(s)\theta(f)\lambda(s)^*)\\
 &= s\cdot \p(f), 
\end{align*}
where we note that 
$E(\lambda(s)x\lambda(s)^*)=s\cdot E(x)$ for any $x\in\vng$. 
It follows from Proposition \ref{prop} that $\p=\id$. 
Since $E$ is faithful, this implies that 
$\theta|_{C(\bG)}=\id_{C(\bG)}$ and the conclusion follows. 
\end{proof}

\begin{proof}[Proof of Corollary \ref{cor}]
Since $\vng$ is injective, the injective envelope $I(\cg)$ 
sits (as an operator system) between $\cg$ and $\vng$. 
By Theorem~\ref{thm}, a completely positive projection 
$\theta$ from $\vng$ onto $I(\cg)$ acts identically on $\cgb$. 
It follows that $\cgb$ is contained in $I(\cg)$. 
This inclusion agrees with each $C^*$-algebra structure, 
because the C$^*$-algebra structure of $I(\cg)$ comes from 
the Choi-Effros product \cite{ce} associated with $\theta$.
\end{proof}
%
%
%

\end{document}